\documentclass[11pt]{article}

\usepackage{amsmath,amsthm}
\usepackage{amssymb}
\usepackage{hyperref}

\allowdisplaybreaks

\begin{document}

\title{Euler sums of generalized hyperharmonic numbers}
\author{
Rusen Li
\\
\small School of Mathematics\\
\small Shandong University\\
\small Jinan 250100 China\\
\small \texttt{limanjiashe@163.com}
}

\date{
\small 2020 MR Subject Classifications: 11B37, 11B68, 11M32
%\small Submitted: November 15, 2019;  Accepted: December 25, 2019.\\
%\small MR Subject Classifications: Primary 11B65; Secondary 05A19
}

\maketitle

\def\stf#1#2{\left[#1\atop#2\right]}
\def\sts#1#2{\left\{#1\atop#2\right\}}
\def\e{\mathfrak e}
\def\f{\mathfrak f}

\newtheorem{theorem}{Theorem}
\newtheorem{Prop}{Proposition}
\newtheorem{Cor}{Corollary}
\newtheorem{Lem}{Lemma}
\newtheorem{Example}{Example}
\newtheorem{Remark}{Remark}
\newtheorem{Definition}{Definition}
\newtheorem{Conjecture}{Conjecture}

\begin{abstract}
In this paper, we mainly show that Euler sums of generalized hyperharmonic numbers can be expressed in terms of linear combinations of the classical Euler sums.
\end{abstract}

\section{Introduction}

The main investigating object of the this paper is the so called Euler sums of generalized hyperharmonic numbers
$$
\zeta_{H^{(p,r)}}(m):= \sum_{n=1}^\infty \frac{H_n^{(p,r)}}{n^{m}}\quad (p, r, m \in \mathbb N:=\{1,2,3,\cdots\})\,,
$$
where
$$
H_n^{(p,r)}:=\sum_{j=1}^n H_{j}^{(p,r-1)} \quad (n, p, r \in \mathbb N)
$$
are the generalized hyperharmonic numbers (see {\cite{Dil,omur}}). Furthermore, $H_n^{(p,1)}=H_n^{(p)}=\sum_{j=1}^{n} 1/n^{p}$ are the generalized harmonic numbers and $H_n^{(1,r)}=h_n^{(r)}$ are the classical hyperharmonic numbers. In particular $H_n^{(1,1)}=H_n$ are the classical harmonic numbers.

Many researchers have been studying Euler sums of harmonic and hyperharmonic numbers (see \cite{Dil,Flajolet,Kamano,mezo} and references therein), since they play important roles in combinatorics, number theory, analysis of algorithms and many other areas (see e.g. \cite{Knuth}). It is interesting that the Riemann zeta function $\zeta(s):=\sum_{n=1}^{\infty}1/n^{s}$ often appears in such expressions. A well-known result \cite{Flajolet} that can be traced back to the time of Euler is as the following:
$$
2\sum_{n=1}^\infty \frac{H_n}{n^{m}}=(m+2)\zeta(m+1)-\sum_{n=1}^{m-2}\zeta(m-n)\zeta(n+1), \quad m=2,3,\cdots.
$$
From that time on, various similar types of infinite series had been investigated. For instance, Flajolet and Salvy \cite{Flajolet} developed the contour integral representation approach to the evaluation of Euler sums involving the classical harmonic numbers. Following Flajolet-Salvy's paper \cite{Flajolet}, we write the classical linear Euler sums as
$$
S_{p,q}^{+,+}:=\sum_{n=1}^\infty \frac{H_n^{(p)}}{{n}^{q}}\,.
$$

Mez\H o and Dil \cite{mezo} considered the infinite sum
$$
\sum_{n=1}^\infty \frac{h_n^{(r)}}{n^{m}}
$$
for an integer $m\ge r+1$, and showed that it could be expressed in terms of infinite sums of the Hurwitz zeta function values, where the well-known Hurwitz zeta function is defined as
$$
\zeta(s,a)=\sum_{n=0}^\infty \frac{1}{(n+a)^{s}} \quad (s\in \mathbb C, \mathfrak Re(s)>1, a>0).
$$
Note that $\mathbb C$ denotes the set of complex numbers and $\mathfrak Re(s)$ denotes the real part of the complex number $s$.

Later Dil and Boyadzhiev \cite{Dil} extended this result to infinite sums involving multiple sums of the Hurwitz zeta function values.

If we regard $\sum_{n=1}^\infty h_n^{(r)}/{n^{s}}$ as a complex function in variable $s$, there are some more progresses toward this direction. For instance, Kamano \cite{Kamano} expressed the complex variable function $\sum_{n=1}^\infty h_n^{(r)}/{n^{s}}$ in terms of the Riemann zeta function, and showed that it could be meromorphically continued to the whole complex plane. In addition, the residue at each pole was also given.

As a natural generalization, Dil, Mez\H o and Cenkci \cite{Dil1} considered Euler sums of generalized hyperharmonic numbers
$$
\zeta_{H^{(p,r)}}(m):= \sum_{n=1}^\infty \frac{H_n^{(p,r)}}{n^{m}}.
$$
They proved that for positive integers $p, r$ and a positive integer $m$ with $m>r$, $\zeta_{H^{(p,r)}}(m)$ could be expressed in terms of series of multiple sums of the Hurwitz zeta function values. For $r=1, 2, 3$, $\zeta_{H^{(p,r)}}(m)$ were also given explicit expressions as linear combinations of the multiple zeta values. On the contrary, \" Om\" ur and Koparal defined two
$n\times n$ matrices $A_n$ and $B_n$ with $a_{i,j}=H_i^{(j,r)}$ and $b_{i,j}=H_i^{(p,j)}$, respectively, and gave some interesting factorizations and determinants of the matrices $A_n$ and $B_n$.

The motivation of this paper arises from Dil, Mez\H o and Cenkci's result (see \cite{Dil1}). Although they reduced $\zeta_{H^{(p,r)}}(m)$ to zeta values for small $p, r, m$, they didn't find a general formula. Our main aim is to establish a general formula to express $\zeta_{H^{(p,r)}}(m)$ in terms of linear combinations of the classical Euler sums. Since a big family of the classical Euler sums can be reduced to zeta values (see \cite{Flajolet} and references therein), we can reduce $\zeta_{H^{(p,r)}}(m)$ to zeta values for appropriate values of $p, r, m$. In addition, we also present several conjectures on these coefficients.

\section{Main theorem}

We are now going to prove our main theorem of this section. Before going further, we introduce some notations and lemmata.

It is well-known that the sum of powers of consecutive intergers $1^k+2^k+\cdots+n^k$ can be explicitly expressed in terms of Bernoulli numbers or Bernoulli polynomials. Faulhaber's formula can be written as
\begin{align}
\sum_{\ell=1}^{n}\ell^{k}&=\frac{1}{k+1}\sum_{j=0}^k \binom{k+1}{j}B_j^{+} n^{k+1-j}\label{ber}\\
                       &=\frac{1}{k+1}(B_{k+1}(n+1)-B_{k+1}(1))\quad\hbox{\cite{CFZ}}\,,
\label{ber1}
\end{align}
where Bernoulli numbers $B_n^{+}$ are determined by the recurrence formula
$$
\sum_{j=0}^k\binom{k+1}{j}B_j^{+}=k+1\quad (k\ge 0)
$$
or by the generating function
$$
\frac{t}{1-e^{-t}}=\sum_{n=0}^\infty B_n^{+}\frac{t^n}{n!}\,,
$$
and  Bernoulli polynomials $B_n(x)$ are defined by the following generating function
$$
\frac{te^{xt}}{e^{t}-1}=\sum_{n=0}^\infty B_n(x)\frac{t^n}{n!}\,.
$$
It is known that the Conway-Guy formula \cite{BGG,Conway} for the hyperharmonic numbers has a simple expression
\begin{align}\label{hypercc}
H_{n}^{(1,r)}=h_n^{(r)}=\binom{n+r-1}{r-1}(H_{n+r-1}-H_{r-1}).
\end{align}
It seems difficult to find a similar formula for the generalized hyperharmonic numbers $H_n^{(p,r)}$, because it looks hard to reduce the $(p, r, n)$ parameter to one parameter as in (\ref{hypercc}), where $r$ was reduced to $r=1$. We make further investigations in this direction, since it is important in the proof of the main theorem. We try to reduce the $(p, r, n)$ parameter to the $(p, n)$ parameter.

\begin{Lem}\label{lemma1}
For $r, n, t \in \mathbb N$, defining
\begin{align}
T(r,n,t):=\sum_{k_{1}=t}^{n} \sum_{k_{2}=t}^{k_{1}} \cdots \sum_{k_{r-1}=t}^{k_{r-2}}1 \,,
\end{align}
then we have
\begin{align}\label{brtm}
T(r,n,t)=\sum_{m=0}^{r-1} B(r,t,m) n^m \,,
\end{align}

where $B(r,t,m)$ satisfy the following recurrence relations
\begin{align}
&B(r+1,t,\ell)=\sum_{m=\ell-1}^{r-1} \frac{B(r,t,m)}{m+1} \binom{m+1}{m-\ell+1}B_{m-\ell+1}^{+} \quad (1\leq \ell \leq r)\,, \label{brt1}\\
&B(r+1,t,0)=-\sum_{j=0}^{r-1}(t-1)^{1+j} \sum_{m=j}^{r-1} \frac{B(r,t,m)}{m+1} \binom{m+1}{m-j}B_{m-j}^{+}\,, \label{brt0}
\end{align}
with boundary value $B(1,t,0)=1$. In addition, $B(r,t,m)$ denote polynomials in variable $t$ of at most $r-m-1$ degree.
\end{Lem}

\begin{proof}
From the definition of $T(r,n,t)$, we have
\begin{align*}
T(r+1,n,t)
&=\sum_{m=0}^{r-1} B(r,t,m) \sum_{k_{1}=t}^{n} k_{1}^m \notag\\
&=\sum_{m=0}^{r-1} B(r,t,m) \frac{1}{m+1}\sum_{j=0}^{m} \binom{m+1}{j}B_j^{+} (n^{m+1-j}-(t-1)^{m+1-j}) \notag\\
&=\sum_{j=0}^{r-1} n^{1+j}\sum_{m=j}^{r-1}B(r,t,m)\frac{1}{m+1}\binom{m+1}{m-j}B_{m-j}^{+} \notag\\
&\quad \quad-\sum_{j=0}^{r-1} (t-1)^{1+j}\sum_{m=j}^{r-1}B(r,t,m) \frac{1}{m+1}\binom{m+1}{m-j}B_{m-j}^{+}\,. \notag\\\
\end{align*}
Since $T(r+1,n,t)=\sum_{m=0}^{r}B(r+1,t,m) n^m$, comparing the coefficients of $n^m$ gives the recurrence relations.

We prove the claim on the degree of $B(r,t,m)$ by induction on $r$.
For $r=1$, from the definition of $T(r,n,t)$, we have $B(1,t,0)=1$. The claim is true for $r=1$. Assume the claim is true for $r$, thus we show that the claim is true for $r+1$.
Since $B(r,t,m)$ is of at most $r-m-1$ degree, with the help of formula (\ref{brt1}), we obtain that for $B(r+1,t,\ell)$ is of at most $r-\ell$ degree. From formulas (\ref{brt1}) and (\ref{brt0}), we have
$$
B(r+1,t,0)=-\sum_{j=0}^{r-1}(t-1)^{1+j}B(r+1,t,j+1).
$$
Since each term is of at most $r$ degree, we get that $B(r+1,t,0)$ is of at most $r$ degree.
\end{proof}

Since $B(r,t,m)$ is of at most $r-m-1$ degree, we can write $B(r,t,m)$ more precisely.
For $r \in \mathbb N$ and $0 \leq m \leq r-1$, let
\begin{align*}
B(r,t,m)=\sum_{j=0}^{r-1-m} b(r,m,j) t^{j} \,.
\end{align*}
Since $B(1,t,0)=1$, we have $b(1,0,0)=1$.
We now give recurrence relations for $b(r,m,j)$ more precisely.

\begin{Lem}\label{lem2}
For $r \in \mathbb N$, one has
\begin{align*}
&b(r+1,\ell,j)=\sum_{m=\ell-1}^{r-1-j} \frac{b(r,m,j)}{m+1} \binom{m+1}{m-\ell+1}B_{m-\ell+1}^{+}\quad(1\leq \ell \leq r, 0\leq j \leq r-\ell)\,,\\
&b(r+1,0,p)=-\sum_{j=0}^{r-1} \sum_{\ell=max\{0, p+1+j-r\}}^{min\{1+j, p\}}C(r,p,j,\ell)\quad (0\leq p \leq r)\,,
\end{align*}
where
$$
C(r,p,j,\ell)=\binom{1+j}{\ell}(-1)^{1+j-\ell} \sum_{m=j}^{r-1-p+\ell} \frac{b(r,m,p-\ell)}{m+1}\binom{m+1}{m-j}B_{m-j}^{+}.
$$

\end{Lem}
\begin{proof}
From the formula (\ref{brt1}), we have
\begin{align*}
B(r+1,t,\ell)
&=\sum_{m=\ell-1}^{r-1} \frac{B(r,t,m)}{m+1} \binom{m+1}{m-\ell+1}B_{m-\ell+1}^{+}\notag\\
&=\sum_{m=\ell-1}^{r-1} \sum_{j=0}^{r-1-m} b(r,m,j)t^{j}\frac{1}{m+1} \binom{m+1}{m-\ell+1}B_{m-\ell+1}^{+}\notag\\
&=\sum_{j=0}^{r-\ell} t^{j} \sum_{m=\ell-1}^{r-1-j}b(r,m,j)\frac{1}{m+1} \binom{m+1}{m-\ell+1}B_{m-\ell+1}^{+}.\notag\\
\end{align*}
Since $B(r+1,t,\ell)=\sum_{j=0}^{r-\ell} b(r+1,\ell,j) t^{j}$, comparing the coefficients of $t^{j}$ gives the first recurrence relation.

From the formula (\ref{brt0}), we have
\begin{align*}
&\quad B(r+1,t,0)\notag\\
&=-\sum_{j=0}^{r-1}(t-1)^{1+j} \sum_{m=j}^{r-1} \frac{B(r,t,m)}{m+1} \binom{m+1}{m-j}B_{m-j}^{+}\notag\\
&=-\sum_{j=0}^{r-1}(t-1)^{1+j} \sum_{m=j}^{r-1} \sum_{\ell=0}^{r-1-m}b(r,m,\ell)t^{\ell}\frac{1}{m+1} \binom{m+1}{m-j}B_{m-j}^{+}\notag\\
&=-\sum_{j=0}^{r-1}\sum_{k=0}^{1+j}\binom{1+j}{k}(-1)^{1+j-k}t^{k} \sum_{\ell=0}^{r-1-j}t^{\ell}\sum_{m=j}^{r-1-\ell} b(r,m,\ell)\frac{1}{m+1} \binom{m+1}{m-j}B_{m-j}^{+}\notag\\
&=-\sum_{p=0}^{r}t^{p} \sum_{j=0}^{r-1} \sum_{\ell=max\{0, p+1+j-r\}}^{min\{1+j, p\}}\binom{1+j}{\ell}(-1)^{1+j-\ell}\notag\\
 &\quad\quad \times \sum_{m=j}^{r-1-p+\ell} \frac{b(r,m,p-\ell)}{m+1}\binom{m+1}{m-j}B_{m-j}^{+}.\notag\\
\end{align*}
Since $B(r+1,t,0)=\sum_{p=0}^{r} b(r+1,0,p) t^{p}$, comparing the coefficients of $t^{p}$ gives the second recurrence relation.
\end{proof}

If we write $a(r,m,j):=b(r,j,m)$, then we have
\begin{align}\label{arnm}
T(r,n,t)=\sum_{m=0}^{r-1} \sum_{j=0}^{r-1-m} a(r,m,j) n^{j} t^{m}.
\end{align}
We give a simple expression for $H_n^{(p,r)}$:
\begin{align}\label{hnpr}
H_n^{(p,r)}
&=\sum_{t=1}^{n}\frac{1}{t^p}\sum_{k_{1}=t}^{n} \sum_{k_{2}=t}^{k_{1}}\cdots \sum_{k_{r-1}=t}^{k_{r-2}}1\,\notag\\
&=\sum_{t=1}^{n}\frac{1}{t^p}T(r,n,t)\,\notag\\
&=\sum_{m=0}^{r-1} \sum_{j=0}^{r-1-m} a(r,m,j) n^{j} H_{n}^{(p-m)}\,.
\end{align}
Note that, when $\ell \geq 0$, $ H_{n}^{(-\ell)}$ is understood to be the sum $\sum_{x=1}^{n} x^{\ell}$.

Now we are able to prove our main theorem.
\begin{theorem}\label{maintheorem}
Let $r, p, m \in \mathbb N$ with $m \geq r+1$, we have,
\begin{align}
\zeta_{H^{(p,r)}}(m)=\sum_{\ell=0}^{r-1} \sum_{j=0}^{r-1-\ell} a(r,\ell,j)S_{p-\ell,m-j}^{+,+}.
\end{align}
Therefore $\zeta_{H^{(p,r)}}(m)$ can be expressed in terms of linear combinations of the classical Euler sums. Moreover, it can be expressed as linear combinations of multiple zeta sums.
\end{theorem}
\begin{proof}
Using formual (\ref{hnpr}), we can write
\begin{align*}
\zeta_{H^{(p,r)}}(m)&=\sum_{n=1}^\infty \frac{1}{n^{m}}\sum_{\ell=0}^{r-1} \sum_{j=0}^{r-1-\ell} a(r,\ell,j) n^{j}H_{n}^{(p-\ell)}\,\\
&=\sum_{\ell=0}^{r-1} \sum_{j=0}^{r-1-\ell} a(r,\ell,j) \sum_{n=1}^\infty \frac{H_{n}^{(p-\ell)}}{n^{m-j}}\,.
\end{align*}
Note that
$$
\sum_{n=1}^\infty \frac{H_{n}^{(p-\ell,1)}}{n^{m-j}}
=\sum_{n=1}^\infty \frac{H_{n-1}^{(p-\ell,1)}}{n^{m-j}}+\zeta(m-j+p-\ell),
$$
and $\sum_{n=1}^\infty H_{n-1}^{(p-\ell,1)}/{n^{m-j}}$ is a specialization of a multiple zeta function (see \cite{Dil1} and references therein), we get the desired result.
\end{proof}

\section{More results on $a(r,m,\ell)$}

Lemma \ref{lem2} gives a recurrence relation for $b(r,m,\ell)$. In this section, We will establish another recurrence relation for $a(r,m,\ell)$.

\begin{theorem}\label{lem3}
For $r \in \mathbb N$, one has
\begin{align*}
&a(r+1,r,0)=-\sum_{m=0}^{r-1} a(r,m,r-m-1)\frac{1}{r-m}\,,\\
&a(r+1,m,\ell)=\sum_{j=\ell-1}^{r-1-m} \frac{a(r,m,j)}{j+1} \binom{j+1}{j-\ell+1}B_{j-\ell+1}^{+}\quad(0\leq m \leq r-1, 1\leq \ell \leq r-m)\,,\\
&a(r+1,m,0)=-\sum_{y=0}^{m} \sum_{j=max\{0, m-y-1\}}^{r-1-y}a(r,y,j)D(r,m,j,y)\quad (0\leq m \leq r-1)\,,
\end{align*}
where
$$
D(r,m,j,y)=\sum_{\ell=max\{0, m-y-1\}}^{j} \frac{1}{j+1} \binom{j+1}{j-\ell}B_{j-\ell}^{+}\binom{\ell+1}{m-y}(-1)^{1+\ell-m+y}.
$$

\end{theorem}
\begin{proof}
From the formula (\ref{arnm}), we have
\begin{align*}
&\quad T(r+1,n,t)\notag\\
&=\sum_{m=0}^{r-1} \sum_{j=0}^{r-1-m} a(r,m,j) t^{m} \sum_{k_{1}=t}^{n} k_{1}^{j}\notag\\
&=\sum_{m=0}^{r-1} \sum_{j=0}^{r-1-m} a(r,m,j) t^{m} \frac{1}{j+1}\sum_{\ell=0}^{j} \binom{j+1}{\ell}B_{\ell}^{+} (n^{j+1-\ell}-(t-1)^{j+1-\ell})\notag\\
&=\sum_{m=0}^{r-1} t^{m} \sum_{j=0}^{r-1-m} a(r,m,j) \frac{1}{j+1}\sum_{\ell=0}^{j} \binom{j+1}{j-\ell}B_{j-\ell}^{+} n^{1+\ell}\notag\\
&\quad -\sum_{m=0}^{r-1} \sum_{j=0}^{r-1-m} a(r,m,j) t^{m} \frac{1}{j+1}\sum_{\ell=0}^{j} \binom{j+1}{j-\ell}B_{j-\ell}^{+} \sum_{x=0}^{1+\ell} \binom{\ell+1}{x}t^{x}(-1)^{1+\ell-x}\notag\\
&=\sum_{m=0}^{r-1} t^{m} \sum_{j=0}^{r-1-m} a(r,m,j) \frac{1}{j+1}\sum_{\ell=0}^{j} \binom{j+1}{j-\ell}B_{j-\ell}^{+} n^{1+\ell}\notag\\
&\quad -\sum_{m=0}^{r-1} \sum_{j=0}^{r-1-m}\frac{a(r,m,j)}{j+1}\sum_{x=0}^{1+j}t^{m+x} \sum_{\ell=max\{0,x-1\}}^{j} \binom{j+1}{j-\ell}B_{j-\ell}^{+} \binom{\ell+1}{x}(-1)^{1+\ell-x}\notag\\
&=\sum_{m=0}^{r-1} t^{m} \sum_{j=0}^{r-1-m} a(r,m,j) \frac{1}{j+1}\sum_{\ell=0}^{j} \binom{j+1}{j-\ell}B_{j-\ell}^{+} n^{1+\ell}\notag\\
&\quad -\sum_{m=0}^{r} t^{m}\sum_{y=0}^{min\{r-1,m\}} \sum_{j=max\{0, m-y-1\}}^{r-1-y}a(r,y,j)D(r,m,j,y).\notag\\
\end{align*}
On the other hand,
$$
T(r+1,n,t)=\sum_{m=0}^{r} \sum_{j=0}^{r-m} a(r+1,m,j) n^{j} t^{m},
$$
comparing the coefficients gives the desired result.
\end{proof}

Next we give explicit values of $a(r,m,j)$ for small $r$.

{\bf Case} $r=1$:
\begin{align*}
&a(1,0,0)=1.
\end{align*}

{\bf Case} $r=2$:
\begin{align*}
&a(2,1,0)=-1,\\
&a(2,0,0)=1, \quad a(2,0,1)=1.
\end{align*}

{\bf Case} $r=3$:
\begin{align*}
&a(3,2,0)=\frac{1}{2},\\
&a(3,1,0)=-\frac{3}{2}, \quad a(3,1,1)=-1,\\
&a(3,0,0)=1, \quad a(3,0,1)=\frac{3}{2}, \quad a(3,0,2)=\frac{1}{2}.\\
\end{align*}

{\bf Case} $r=4$:
\begin{align*}
&a(4,3,0)=-\frac{1}{6},\\
&a(4,2,0)=1, \quad a(4,2,1)=\frac{1}{2},\\
&a(4,1,0)=-\frac{11}{6}, \quad a(4,1,1)=-2, \quad a(4,1,2)=-\frac{1}{2},\\
&a(4,0,0)=1, \quad a(4,0,1)=\frac{11}{6}, \quad a(4,0,2)=1, \quad a(4,0,3)=\frac{1}{6}.\\
\end{align*}

{\bf Case} $r=5$:
\begin{align*}
&a(5,4,0)=\frac{1}{24},\\
&a(5,3,0)=-\frac{5}{12}, \quad a(5,3,1)=-\frac{1}{6},\\
&a(5,2,0)=\frac{35}{24}, \quad a(5,2,1)=\frac{5}{4}, \quad a(5,2,2)=\frac{1}{4},\\
&a(5,1,0)=-\frac{25}{12}, \quad a(5,1,1)=-\frac{35}{12}, \quad a(5,1,2)=-\frac{5}{4}, \quad a(5,1,3)=-\frac{1}{6},\\
&a(5,0,0)=1,\quad a(5,0,1)=\frac{25}{12}, \quad a(5,0,2)=\frac{35}{24}, \quad a(5,0,3)=\frac{5}{12}, \quad a(5,0,4)=\frac{1}{24}.\\
\end{align*}

Observing the above facts, we present the following four conjectures.
\begin{Conjecture}
For $r, m, \ell \in \mathbb N$ with $0\leq m \leq r-1$ and $0\leq \ell \leq r-1-m$, we conjecture that
$$
a(r,m,\ell)=(-1)^{m+\ell} a(r,\ell,m)\,.
$$
\end{Conjecture}

\begin{Conjecture}
For $r \in \mathbb N$ with $0\leq n \leq r-1$, we conjecture that
$$
\sum_{\ell=0}^{n}a(r,n-\ell,\ell)=\delta_{n 0},
$$
where $\delta_{n m}$ is the Kronecker delta, that is, $\delta_{n n}=1$, $\delta_{n m}=0$ for $n \neq m$.
\end{Conjecture}

\begin{Conjecture}
For $r \in \mathbb N$, we conjecture that
$$
\sum_{\ell=0}^{r-1}a(r,0,\ell)=r \quad (r \geq 1),
$$
and
$$
\sum_{\ell=0}^{r-1}a(r,\ell,0)=0 \quad (r \geq 2).
$$
\end{Conjecture}

\begin{Conjecture}
For $r, m \in \mathbb N$ with $0\leq m \leq r-1$, we conjecture that
$$
\mathrm{sgn}(a(r,m,\ell))=(-1)^{m} \quad (0\leq \ell \leq r-1-m),
$$
where $\mathrm{sgn}(x)$ is the signum function defined by
\begin{equation*}
\mathrm{sgn}(x)=
\begin{cases}
1& \text{$x>0$},\\
0& \text{$x=0$},\\
-1& \text{$x<0$}.
\end{cases}
\end{equation*}
In particular
$$
a(r,m,\ell)\neq 0 \quad (0\leq \ell \leq r-1-m).
$$
\end{Conjecture}

\end{document}